\documentclass[a4paper,11pt,twoside]{article}

\usepackage{xypic}
\usepackage{latexsym}
\usepackage{amsmath}
\usepackage{amstext}
\usepackage{amsfonts}
\usepackage{amssymb}
\usepackage{amsbsy}
\usepackage{enumerate}
\usepackage{theorem}

\xyoption{all}
\CompileMatrices

\newtheorem{lemma}{Lemma}
\newtheorem{proposition}[lemma]{Proposition}

\newtheorem{notation}[lemma]{Notation}

\theorembodyfont{\rmfamily}

\newtheorem{definition}[lemma]{Definition}

\newcommand{\ns}{\mathop{\rm ns}\nolimits}
\newcommand{\s}{\mathop{\rm s}\nolimits}

\newcommand{\Div}{\mathop{\rm Div}\nolimits}
\newcommand{\Cl}{\mathop{\rm Cl}\nolimits}

\begin{document}

\addtolength{\baselineskip}{0.5 pt}

\setlength{\parskip}{0.2 ex}

\title{Index calculus with double large prime variation for curves of small genus with cyclic class group}

\author{Claus Diem}

\date{June 23, 2006}

\maketitle

\begin{abstract}
We present an index calculus algorithm with double large prime variation which lends itself well to a rigorous analysis. Using this algorithm we prove that for fixed genus $g \geq 2$, the discrete logarithm problem in degree 0 class groups of non-singular curves over finite fields $\mathbb{F}_q$ can be solved in an expected time of $\tilde{O}(q^{2-2/g})$, provided that the curve is given by a plane model of bounded degree and the degree 0 class group is cyclic.

The result generalizes a previous result for hyperelliptic curves given by an imaginary Weierstra{\ss} equation obtained by Gaudry, Thom\'e, Th\'eriault and the author.
\end{abstract}


\begin{center}\begin{minipage}{110mm}\footnotesize{\bf MSC2000:} Primary:  11Y16; Secondary: 14G50, 11G20

\end{minipage}
\end{center}
\vspace{1 ex}

\pagestyle{myheadings} \markboth{\sc Diem}{\sc Index calculus with double large prime variation}

\section{Introduction}\label{intro}

In \cite{GTTD} and \cite{Nag} index calculus algorithms with double large prime variation for the solution of the discrete logarithm problem (DLP) in degree 0 class groups of curves of small genus have been given. In \cite{GTTD} it has been proven that for fixed genus $g \geq 2$ one can solve the DLP in degree 0 class groups of hyperelliptic curves over finite fields $\mathbb{F}_q$, given by an imaginary Weierstra{\ss} equation, in an expected time of $\tilde{O}(q^{2-2/g})$, provided that the degree 0 class group is cyclic or its structure is known.

In this work, we generalize this result from hyperelliptic curves given by an imaginary Weierstra{\ss} equation to arbitrary curves. We thereby keep the restriction that the degree 0 class group is cyclic or its structure is known. Our result is as follows.

\paragraph{Theorem}
\emph{Let a natural number $g \geq 2$ be fixed. Then the discrete logarithm problem in cyclic degree 0 class groups of curves of genus $g$ in degree 0 class groups of curves given by plane models of bounded degree can with a randomized algorithm be solved in an expected time of $\tilde{O}(q^{2-2/g})$.}

\smallskip

The $\tilde{O}$-notation captures logarithmic factors. The result also holds if the degree 0 class groups is not necessarily cyclic but its structure is known. By this we mean that a basis of the degree 0 class group as well as the orders of the basis elements are known. If the algorithm is applied with respect two elements of $\Cl^0(\mathcal{C})$ for which the discrete logarithm problem is unsolvable, it outputs ``unsolvable'' in an expected running time of $\tilde{O}(q^{2-2/g})$.

The algorithm has a storage requirement of $\tilde{O}(q^{1- \frac{1}{g} + \frac{1}{g^2}})$. More concretely, although the algorithm is randomized, there exists a function in $\tilde{O}(q^{1-\frac{1}{g} + \frac{1}{g^2}})$ such that the storage requirements are bounded by this function for every run.

In \cite[Theorem 56]{Hess-subexp} it is shown that there exists a constant $C \in \mathbb{N}$ such that any curve of genus $g$ over any finite field can be represented by a plane model of degree $\leq C \cdot g$. This shows that the theorem in fact applies to all curves of a fixed genus whose degree 0 class group is cyclic, provided that the curve is represented appropriately.

\smallskip
The present work is motivated by the fact that the discrete logarithm problem in degree 0 class groups of curves is a well established cryptographic primitive.

We note that the scope of the present work lies purely in the realm of theoretical cryptology / theoretical computational mathematics. For practical computations in degree 0 class groups of non-hyperelliptic curves of small genus, we advice to try to construct a plane model of degree $g+1$ or smaller (see \cite[Section 6]{Di-sd}) and then to follow the algorithm outlined in \cite[Section 5]{Di-sd}.

\smallskip
The algorithm is given in the next section, and the analysis is given in Section \ref{sec-analysis}.

\subsubsection*{Terminology and data structures}

In the theorem, we implicitly used the following terminology and the following conventions concerning data structures.

If not stated otherwise, a \emph{curve} is always irreducible and non-singular. A \emph{plane model} of a curve over a field $k$ is a possibly singular curve in $\mathbb{P}^2_k := \text{Proj}(k[X,Y,Z])$ which is birational to the curve. We represent a curve via a (fixed) plane model, and we represent the plane model (and thus the curve itself) by a defining homogeneous polynomial $F(X,Y,Z)$.

We represent the points of $\mathcal{C}$ by their corresponding points of the plane model, with some additional information for the singular points.

By a \emph{divisor} on a curve $\mathcal{C}$ over a field $k$ we always mean a divisor of $\mathcal{C}$ over $k$ (a $k$-rational divisor). We think of divisors as being represented as a formal sum of closed points of $\mathcal{C}$. (This is called the \emph{free representation} in~\cite{Hess-RR}.) If $D$ is a divisor, the corresponding divisor class is denoted by $[D]$. The degree 0 class group of $\mathcal{C}$ over $\mathbb{F}_q$ is denoted by $\Cl^0(\mathcal{C})$.

Let us consider curves over finite fields $\mathbb{F}_q$. For fixed genus and $q \gg 0$, for any curve $\mathcal{C}$ over $\mathbb{F}_q$, $\mathcal{C}(\mathbb{F}_q)$ is non-empty. We assume that this is the case and fix a point $P_0 \in \mathcal{C}(\mathbb{F}_q)$. After having fixed $P_0$, the elements in $\Cl^0(\mathcal{C})$ can be represented in the following way: An effective divisor $D$ is called \emph{maximally reduced along} $P_0$ if the linear system $|D-P_0|$ is empty. By the Riemann-Roch theorem, maximally reduced divisors have degree $\leq g$, and $D \mapsto [D] - \deg(D) \cdot [P_0]$ defines a bijection between the set of maximally reduced effective divisors and $\Cl^0(\mathcal{C})$ (see \cite[Prop.~8.2.]{Hess-RR}). We think of degree 0 divisor classes  as being represented by their corresponding maximally reduced effective divisors.

With this representation of the elements of the degree 0 class group, the arithmetic in $\Cl^0(\mathcal{C})$ can -- for curves represented by plane models of bounded degree -- be carried out in randomized polynomial time (cf.\ e.g.\ \cite{Vo}, \cite{HI-RR}, \cite{Hess-RR}, \cite{Mak}, \cite{Mak2}).

\subsubsection*{Double large prime variation and this work}

As has already been stated, the theorem has already been proven for hyperelliptic curves given by an imaginary Weierstra{\ss} equation in \cite{GTTD}. Moreover it has been pointed out in \cite{GTTD} that heuristically the result also holds for arbitrary curves of represented by plane models of bounded degree. Because of the close relationship between this work and \cite{GTTD}, we advice the reader to have \cite{GTTD} at hand when he goes through the details of the algorithm.

Let us recall some basic ideas about index calculus with double large prime variation and the proof of the theorem for hyperelliptic curves given by an imaginary Weierstra{\ss} equation in \cite{GTTD}.

Generally speaking, a \emph{double large prime variation} of an index calculus algorithm consists of the following: One not only considers relations which split over the factor base but also takes relations with up to two \emph{large primes} into account. These relations are stored as edges in a so-called \emph{graph of large prime relations}. This graph is used to obtain ``recombined'' relations over the factor base.

There are two double large prime variation algorithms presented in \cite{GTTD}: a ``full algorithm'' and a ``simplified algorithm''. The theorem in \cite{GTTD} is proven with the simplified one. Here one does not construct the whole graph of large prime relations but only a tree; we call such a tree a \emph{tree of large prime relations}. At a later stage any relation which splits into elements of the factor base or vertices of the tree is used to obtain a relation over the factor base.

The main challenge resides in controlling the growth of the tree of large prime relations as well as its depth, that is, the maximal distance of any vertex to the root. For hyperelliptic curves given by an imaginary Weierstra{\ss} equation this is relatively easy because one has a concrete description of the effective divisors which are maximally reduced along the point at infinity, and one therefore knows that the growth process is very regular.

\medskip

The algorithm given in this work is a modification of the ``simplified algorithm'' in \cite{GTTD}. The analysis relies on the following proposition.

\begin{proposition}
\label{number-special-divisors}
For curves of fixed genus $g$ over finite fields $\mathbb{F}_q$, the number of special effective divisors of degree $g$ is in $O(q^{g-1})$.
\end{proposition}
Recall that an effective divisor $D$ is called \emph{special} if the linear system $|K-D|$ is non-empty, where $K$ is a canonical divisor. Note that by the Riemann-Roch theorem, an effective divisor of degree $g$ is non-special if and only if it is the only if the linear system $|D|$ merely contains $D$ itself.

We have a canonical injection from the set of non-special divisors of degree $g$ into the set of along $P_0$ maximally reduced effective divisors: Let $D$ be a non-special effective divisor, and let $D'$ be the unique effective divisor of minimal degree with $D' + (\deg(D) - \deg(D')) \cdot P_0 = D$. Then $D'$ is maximally reduced along $P_0$.

We assume that Proposition \ref{number-special-divisors} is well known to many experts in curves and function fields. For the lack of a suitable reference we give a proof in an appendix. Note that a straightforward application of the Hasse-Weil Bound merely gives that the number in question is in $O(q^{q-{1/2}})$.

This proposition makes it possible to discard all special divisors in the analysis of the construction of the tree of large prime relations. It remains the problem to control the growth of the tree. For this, we modify the ``simplified algorithm'' in \cite{GTTD} in such a way that the depth of the tree (not only the expected value of the depth) always lies in $O(\log(q))$.

\section{The algorithm}

Let $g \in \mathbb{N}$, $g \geq 2$ be fixed.

In the following, we consider the discrete logarithm problem in cyclic degree 0 class groups of curves of genus $g$ over finite fields $\mathbb{F}_q$, given by plane models of bounded degree. We thereby implicitly use the data structures described in the introduction. As stated in the introduction, the theorem also holds if the degree 0 class group is not necessarily cyclic but its structure is known. In this case the algorithm should be modified according to the description in \cite[Section 7]{EG}.

The $L$-polynomial of a curve over $\mathbb{F}_q$ given by a plane model of bounded degree can be computed in (deterministic) polynomial time in $\log(q)$. (This follows from \cite[Theorem H]{Pi-Curves} which in turn relies on Pila's extension of the point counting algorithm by Schoof (\cite{Schoof-EC}) to abelian varieties (\cite{Pi-AV}).) This means in particular that the order of the degree 0 class group can be computed in polynomial time in $\log(q)$. Moreover, the order can be factored in subexponential time with the algorithm in \cite{LP}. In the following, we therefore assume that the order and its factorization are known.

The algorithm below terminates in a finite expected time and solves -- if possible -- the discrete logarithm problem for $q \gg 0$. In order to obtain an algorithm which always yields the solution to the DLP (or outputs ``unsolvable''), one can let the algorithm run ``in parallel'' with a brute force calculation (that is, for every step of the algorithm below, one brute force try is performed).

\smallskip

Let $\mathcal{C}$ be a curve of genus $g$ over $\mathbb{F}_q$ such that $\Cl^0(\mathcal{C})$ is cyclic, and let $a, b \in \Cl^0(\mathcal{C})$. The goal is to determine if $b \in \langle a \rangle$ and -- if this is the case -- to compute an $x \in \mathbb{N}$ such that $x \cdot a = b$.

\smallskip

Let $\ell := \# \Cl^0(\mathcal{C})$.
\subsubsection*{Reduction to the DLP with respect to a generator}

The first step of the algorithm consists of a reduction of the problem to the discrete logarithm problem with respect to a generator of the group:

By \cite[Theorem 34]{Hess-RR}, for $q \gg 0$, there exists some $P \in \mathcal{C}(\mathbb{F}_q) - \{ P_0 \}$ such that $c := [P] - [P_0]$ generates $\Cl^0(\mathcal{C})$.

To find such a point $P$, we iterate over all elements of $\mathcal{C}(\mathbb{F}_q)$. For each $P$, we test for each prime factor $\mu$ of $\ell$ if $\frac{\ell}{\mu} \cdot c \neq0$. If this is the case, the order of $c$ is $\ell$, and we fix the point $P$.

If we have found such a $P$, we proceed as follows: We determine $x_a, x_b \in \mathbb{Z}/\ell\mathbb{Z}$ with $x_a \cdot c = a, x_b \cdot c = b$ with the algorithm described below. Then we try to determine an $x \in \mathbb{Z}/\ell\mathbb{Z}$ with $x \cdot x_a = x_b$. If no such $x$ exists, we output ``unsolvable'', otherwise, we output $x$.

If no such $P$ exists, we try to solve the DLP with brute force.

\smallskip

From now on, we assume that $a$ generates $\Cl^0(\mathcal{C})$.
\subsubsection*{The factor base}
We fix any \emph{factor base} $\mathcal{F} = \{ F_1, F_2, \ldots \} \subset \mathcal{C}(\mathbb{F}_q) - \{ P_0 \}$ of size $\lceil  q^{1-\frac{1}{g-1}} \rceil$. (If $\mathcal{C}(\mathbb{F}_q) - \{ P_0 \}$ contains less elements, the algorithm terminates.) As in \cite{GTTD} let $\mathcal{L} := \mathcal{C}(\mathbb{F}_q) - ( \mathcal{F} \cup \{ P_0 \})$ be the set of \emph{large primes}.
\subsubsection*{Construction of the tree of large prime relations}
Similarly to the ``simplified algorithm'' in \cite{GTTD}, we construct a tree of large prime relations on $\mathcal{L} \stackrel{_\cdot}{\cup} \{ * \}$.

For this we repeatedly choose uniformly randomly $\alpha, \beta \in \mathbb{Z}/\ell \mathbb{Z}$ and compute the along $P_0$ maximally reduced effective divisor $D$ with
\begin{equation}
\label{relation}
[D] - \deg(D) \cdot [P_0] = \alpha a + \beta b \; .
\end{equation}
We thereby choose $\alpha$ and $\beta$ independently of each other and independently of all previous choices.

\smallskip

Recall the following usual definition for double large prime variation algorithms (cf.\ \cite{GTTD}).

\begin{definition}
A relation as (\ref{relation}) is called \emph{Full} if $D$ splits into divisors of the factor base. It is called \emph{FP} if $D$ is completely split and is the sum of elements of the factor base and the non-zero multiple of one large prime. It is called \emph{PP} if $D$ is completely split and is the sum of elements of the factor base and non-zero multiples of two large primes.
\end{definition}

FP relations are stored in the tree of large prime relations by inserting a labeled edge between $*$ and the large prime in the relation, and PP relations are stored by inserting a labeled edge between the two large primes in the relations.

In comparison to the ``simplified algorithm'' in \cite{GTTD}, we modify the construction of the tree of large prime relations:

We construct the tree in \emph{stages}, and during each stage we only attach edges to the tree which are connected to vertices constructed \emph{in the previous stage}. In Stage 1, we attach $\lceil q^{1 -1/g} \rceil$ edges coming from FP relations to the root $*$. Thereafter, we terminate Stage $s$ and start Stage $s+1$ whenever the tree has $2^{s-1} \cdot \lceil q^{1 -1/g} \rceil$ edges.

\smallskip

Let us fix this notation.

\begin{notation}
The set of vertices of a tree $T$ is also denoted by $T$.
\end{notation}

A (semi-)formal description of the construction of the tree is as follows.
\medskip

{\sf
\subsubsection*{Algorithm: Construction of the tree of large prime relations}

Construct a tree on $\mathcal{L} \stackrel{_\cdot}{\cup} \{ * \}$ as follows:\\
Let $T_0$ consist only of the root $*$.\\
Let $s \longleftarrow 1$.\\
Repeat\\
\text{ } \hspace{3 ex} Construct a tree $T_s$ which contains $T_{s-1}$ as a subtree as follows:\\
\text{ } \hspace{3 ex} Repeat\\
\text{ } \hspace{7 ex} Choose $\alpha, \beta \in \mathbb{Z}/\ell \mathbb{Z}$ uniformly and independently \\
\text{ } \hspace{7 ex} of each other and of all previous choices at random.\\
\text{ } \hspace{7 ex} Compute the along $P_0$ maximally reduced divisor $D$ with \\ 
\text{ } \hspace{7 ex} $[D] - \deg(D) \cdot [P_0] = \alpha a + \beta b$.\\
\text{ } \hspace{7 ex} If $D$ splits as $D = \sum_j r_j F_j + c_P P + c_Q Q$ where $c_P > 0, c_Q > 0$, \\ 
\text{ } \hspace{7 ex} $P \in  \mathcal{F} \cup T_{s-1}$ and $Q \in \mathcal{L} - (\mathcal{F} \cup T_s)$, \\
\text{ } \hspace{10 ex} if $P \in \mathcal{F}$ (i.e. if we have an FP relation), \\
\text{ } \hspace{13 ex} insert an edge from $*$ to $Q$ into the tree $T_s$,\\
\text{ } \hspace{10 ex} if $P \in T_{s-1}$ (i.e. if we have a PP relation), \\
\text{ } \hspace{13 ex} insert an edge from $P$ to $Q$ into the tree $T_s$, \\
\text{ } \hspace{10 ex} in both cases labeled with $(r_j)_j$ (in sparse representation).\\
\text{ } \hspace{3 ex}  Until $T_s$ contains $2^{s-1} \cdot \lceil q^{1 - 1/g} \rceil$ edges.\\
\text{ } \hspace{3 ex} Let $s \longleftarrow s +1$.
}

\bigskip

This construction of the tree guarantees that the depth of the tree is always in $O(\log(q))$ (see also inequality (\ref{s-bound}) in the next section).
The main difficulty of the analysis of the algorithm resides in proving that a  tree of sufficient size can be constructed in an expected time of $\tilde{O}(q^{2-2/g})$. This is verified in the next section.

The construction of the tree is abandoned if a predefined number of edges $N_{\max}$ is reached. (To improve the readability we did not include this criterion in the description above.) We could for example set $N_{\max} := \lceil q/4 \rceil$. We will however argue in the analysis of the algorithm in the next section that $N_{\max} := \lceil q^{1 - 1/g + 1/g^2} \rceil$ suffices. This smaller value of $N_{\max}$ only lowers the time for the construction of the tree by a constant factor but decreases the storage requirements substantially. This is analogous to the situation in~\cite{GTTD}.

\smallskip

After a tree $T := T_s$ with $N_{\max}$ edges has been constructed, we proceed as in Phase~2 of the ``simplified algorithm'' in \cite{GTTD}. Let us for simplicity assume that $\ell$ is prime. In the general case (in particular if $\ell$ is not square free), the following description should be modified according to the description in Sections 3 and 4 in \cite{EG}.

\subsubsection*{Construction of the matrix}

We construct a sparse matrix with $\# \mathcal{F}$ columns and $\# \mathcal{F} +1$ rows as follows.

We again generate relations (\ref{relation}) by choosing $\alpha$ and $\beta$ independently uniformly at random. If the divisor $D$ in such a relation splits over $\mathcal{F} \cup T$, we use the tree to substitute the large primes involved by sums of possibly negative multiples of elements of the factor base, and we store the coefficient vector of the ``left-hand side'' of the resulting relation as a new row of a sparse matrix $R$.

\subsubsection*{Linear algebra}

After having constructed $R$, the DLP can be solved via a linear algebra computation as usual in index calculus algorithms for cyclic groups: We compute a vector $v$ over $\mathbb{Z}/\ell \mathbb{Z}$ with $v R = 0$ which is uniformly randomly distributed over all vectors in the kernel of $R^t$ with an algorithm from sparse linear algebra. If then $\sum_i v_i \beta_i \in (\mathbb{Z}/\ell\mathbb{Z})^*$, the solution to the DLP is
\[ x :=  - \frac{\sum_i v_i \alpha_i}{\sum_i v_i \beta_i } \in \mathbb{Z}/\ell\mathbb{Z} \; ,\]
as usual. If the condition $\sum_i v_i \beta_i \in (\mathbb{Z}/\ell\mathbb{Z})^* $ is not satisfied, one could compute a new row for the matrix $R$ and then perform the linear algebra computation again. Repeating this procedure would however with a very small probability lead to a storage requirement which is not in $\tilde{O}(q^{1 - \frac{1}{g} + \frac{1}{g^2}})$.

Because we want to bound the storage requirements for every run (instead of merely bounding the expected value of the storage requirements), we do not insert a new row into the matrix if the computation fails but instead restart the whole computation of the matrix $R$. (The same approach has been taken in \cite{EG}.)

\section{Analysis}

\label{sec-analysis}

We now show that the algorithm outlined above computes a solution to the DLP in an expected time of $\tilde{O}(q^{2-2/g})$ (as always for fixed genus $g \geq 2$ and $q \longrightarrow \infty$). 

As already pointed out, the analysis relies on Proposition \ref{number-special-divisors}. Let $C > 0$ be such that for all curves of genus $g$ over any finite fields $\mathbb{F}_q$ the number of special divisors of degree $g$ is $\leq C \cdot q^{g-1}$.

As in the previous section, let $N_{\max} := \lceil q^{1- 1/g + 1/g^2} \rceil$ be the number of edges (that is, the number of vertices different from $*$) at which the construction of the tree is stopped.

The conditions
\[ \begin{array}{ll}
N_{\max} + \# \mathcal{F} \leq q/4 \; & \quad  \quad \quad \# (\mathcal{C}(\mathbb{F}_q) - \{ P_0 \}) \in [\max\{q^{1-\frac{1}{g}}, q/2 \}, 2q] \; \\[1.5 ex]
\# \Cl^0(\mathcal{C}) \leq 2 q^g \; & \quad  \quad \quad q \geq (4 \cdot g! \cdot C)^g
\end{array} \]
hold for $q \gg 0$; we assume that they are satisfied.

\subsubsection*{Construction of the generator and of the factor base}

One can iterate over all points of $\mathcal{C}(\mathbb{F}_q)$ in an expected time of $\tilde{O}(q)$ as follows:

One iterates over the $(X,Z)$-coordinates, and for each iteration one factors the polynomial describing the $Y$-coordinates of the points with the prescribed $(X,Z)$-coordinates.

With this procedure one can find a point $P \in \mathcal{C}(\mathbb{F}_q)$ such that $[P] - [P_0]$ generates $\Cl^0(\mathcal{C})$ as well as determine an appropriate factor base in an expected time of $\tilde{O}(q)$.

\subsubsection*{Construction of the tree of large prime relations}
We come to the analysis of the growth of the tree of large prime relations.

Note first that by our assumption that $\langle a \rangle = \Cl^0(\mathcal{C})$, if $\alpha$ and $\beta$ are drawn uniformly at random from $\mathbb{Z}/\ell\mathbb{Z}$, $\alpha a + \beta b$ is also drawn uniformly at random from $\Cl^0(\mathcal{C})$. This means that the divisor $D$ in (\ref{relation}) is drawn uniformly at random from the set of all effective divisors which are maximally reduced along $P_0$.

By our assumptions on $q$, we always have
\begin{equation}
\label{q/4-bound}
\# (\mathcal{C}(\mathbb{F}_q) - (T_s \cup \mathcal{F} \cup \{ P_0 \} ) \geq q/2 - q/4 = q/4 \; .
\end{equation}

Let $\Div^g(\mathcal{C})$ be the set of effective divisors of degree $g$ on $\mathcal{C}$, and let $\Div^{g,\ns}(\mathcal{C})$ (resp.\ $\Div^{g,\s}(\mathcal{C})$) be the subset of non-special (resp.\ special) effective divisors of degree $g$.

\smallskip

Let us first assume that we are still in Stage 1, that is, only relations with one large prime (not yet in the tree) are considered.

Let us assume we are given the tree $T_0$ or we have already constructed a tree $T_1$ with $< \lceil q^{1 - 1/g} \rceil$ edges. We want to bound the expected number of relations (\ref{relation}) needed until a new edge is inserted into the tree.

Let
\[ \begin{array}{l}
S := \big\{ P_1 + \cdots + P_g \in \Div^{g}(\mathcal{C}) \, | \; \forall i=1, \ldots, g-1 : P_i \in \mathcal{F} , \\[0.5 ex]
\; \; \; \quad \quad \quad \quad \quad \quad \quad \quad \quad \quad \quad \quad \quad P_g \in \mathcal{C}(\mathbb{F}_q) - (T_1 \cup \mathcal{F} \cup \{ P_0 \}) \big\} \, , 
\end{array} \]
\[ S^{\ns} := S \cap \Div^{g,\ns}(\mathcal{C}) \; . \]

Note that any divisor $D \in S^{\ns}$ is maximally reduced along $P_0$ (because $P_0$ is not contained in the support of $D$ and the linear system $|D|$ consists merely of $D$). If a divisor $D = P_1 + \cdots + P_g$ as in the set $S^{\ns}$ appears in a relation~(\ref{relation}), a new edge is inserted into the tree. (Other divisors might also lead to new edges: We ignore FP relations which involve a larger multiple of the large prime, we ignore non-special divisors, and we ignore divisors of degree $< g$.)

We have
\begin{equation}
\label{number-S}
\begin{array}{l}
\# S = \binom{\# \mathcal{F} + g -2}{g-1} \cdot \# \big( \mathcal{C}(\mathbb{F}_q) - (T_1 \cup \mathcal{F} \cup \{ P_0 \}) \big) \geq \frac{\# \mathcal{F}^{g-1}}{(g-1)!} \cdot q/4  \quad \quad \text{ by } (\ref{q/4-bound})\\
\text{ } \quad \; \geq \frac{1}{4 (g-1)!} \cdot q^{\frac{(g-1)^2}{g}} \cdot q = \frac{1}{4 (g-1)!} \cdot q^{\frac{g^2 - g +1}{g}} = \frac{1}{4 (g-1)!} \cdot q^{g - 1 + \frac{1}{g}} \; . 
\end{array} 
\end{equation}

By our assumption that $q \geq (4 \cdot g! \cdot C)^g$, we have
\begin{equation}
\label{number-ns}
\# \Div^{g,\s}(\mathcal{C}) \leq C q^{g-1} \leq \frac{1}{4 g!} \cdot q^{g-1 + \frac{1}{g}} \leq \frac{1}{8 (g-1)!} \cdot q^{g-1 + \frac{1}{g}} \; .
\end{equation}

Inequalities (\ref{number-S}) and (\ref{number-ns}) imply
\[ \# S^{\ns} \geq \frac{1}{8 (g-1)!} \cdot q^{g -1 + \frac{1}{g}} \; . \]
Together with our assumption that $\# \Cl^0(\mathcal{C}) \leq 2q^g$, this implies that the probability that a relation (\ref{relation}) enlarges the tree is 
\[ \geq \frac{\# S^{\ns}}{\# \Cl^0(\mathcal{C})} \geq \frac{1}{16 (g-1)!} \cdot q^{-(1-\frac{1}{g})} \; . \]
The expected number of relations (\ref{relation}) which have to be considered until the tree is enlarged is thus
\[ \leq 16 (g-1)! \cdot q^{1 - \frac{1}{g}} \; .\]
This implies that the expected number of tries until the tree has $\lceil q^{1-\frac{1}{g}} \rceil$ edges is 
\[ \leq 16 (g-1)! \cdot q^{1-\frac{1}{g}} \cdot \lceil q^{1-\frac{1}{g}} \rceil \leq 16 (g-1)! \cdot (q+1)^{2 - \frac{2}{g}} \; . \]

\smallskip

We now assume that $s \geq 2$ and a tree $T_{s-1}$ with $2^{s-2} \cdot \lceil q^{1 - 1/g} \rceil$ edges and a tree $T_s$ with $< 2^{s-1} \cdot \lceil q^{1 - 1/g} \rceil$ edges has already been constructed. The task is again to derive a bound on the expected number of relations (\ref{relation}) needed until the tree is enlarged.

Similarly to above, let
\[ \begin{array}{l}
S := \big\{  P_1 + \cdots + P_g \in \Div^{g}(\mathcal{C}) | \; \forall i=1, \ldots, g-2 : P_i \in \mathcal{F}, \\[0.5 ex]
\text{ } \; \; \; \quad \quad \quad \quad \quad \quad \quad \quad  P_{g-1} \in \mathcal{F} \cup T_{s-1}, \; P_g \in \mathcal{C}(\mathbb{F}_q) - (T_s \cup \mathcal{F} \cup \{ P_0 \}) \big\} \; , 
\end{array} \]
\[ S^{\ns} := S \cap \Div^{g,\ns}(\mathcal{C}) \; . \]
We now have
\begin{equation}
\label{number-S-later}
\begin{array}{l}
\# S = \big( \binom{\# \mathcal{F} + g -2}{g-1} + \binom{\# \mathcal{F} + g -3}{g-2} \cdot \#( T_{s-1} - \{ * \} ) \big) \cdot \# \big( \mathcal{C}(\mathbb{F}_q) - (T_s \cup \mathcal{F} \cup \{ P_0 \}) \big) \\[1 ex]
\text{ } \quad \; \geq ( \frac{\# \mathcal{F}^{g-1}}{(g-1)!} + \frac{\# \mathcal{F}^{g-2}}{(g-2)!} \cdot 2^{s-2} \cdot q^{1 - 1/g})  \cdot q/4  \\[1 ex]
\text{ } \quad \; \geq (\frac{1}{(g-1)!} \cdot q^{\frac{(g-1)^2}{g}} + \frac{1}{(g-2)!} \cdot 2^{s-2} \cdot q^{\frac{(g-1)^2}{g}} ) \cdot q/4 \\[1 ex]
\text{ } \quad \; = ( \frac{1}{4 (g-1)!} + \frac{1}{4 (g-2)!} \cdot 2^{s-2} )\cdot q^{g -1 + \frac{1}{g}} \; .
\end{array} 
\end{equation}
Together with (\ref{number-ns}), this implies
\[ \# S^{\ns} \geq \frac{1}{4 (g-2)!} \cdot 2^{s-2} \cdot q^{g -1 + \frac{1}{g}} \; .\]
This implies that the probability that a relation (\ref{relation}) enlarges the tree is
\[ \geq \frac{1}{8 (g-2)!} \cdot 2^{s-2} \cdot q^{-(1 - \frac{1}{g})} \; . \]
The expected number of relations (\ref{relation}) which have to be considered until the tree is enlarged is thus
\[ \leq 8 (g-2)! \cdot \frac{1}{2^{s-2}} \cdot q^{1 - \frac{1}{g}} \; . \]
This implies that given any tree $T_{s-1}$ with $2^{s-2} \cdot \lceil q^{1 - \frac{1}{g}} \rceil$ edges, the expected number of tries until a tree $T_s$ with $2^{s-1} \cdot \lceil q^{1 - \frac{1}{g}} \rceil$ edges is constructed is 
\[ \leq 16 (g-2)! \cdot (q + 1)^{2 - \frac{2}{g}} \; . \]

\smallskip

We always have have $s = O(\log(q))$ as can be easily be seen: For every run of the algorithm we have for $s \geq 2$
\[ 2q \geq \# (T_s - \{ * \}) \geq \# (T_{s-1} - \{ * \}) = 2^{s-2} \cdot \lceil q^{1 - \frac{1}{g}} \rceil \; ,\]
i.e.
\begin{equation}
\label{s-bound}
s \leq \log_2(q^{\frac{1}{g}}) + 3 = \frac{1}{\log(2) \cdot g} \cdot \log(q) + 3 \; = O(\log(q)) \; .
\end{equation}

It follows that in total an expected number of $O(\log(q) \cdot q^{2 - \frac{2}{g}})$ relations (\ref{relation}) have to be considered until the tree has $N_{\max}$ edges. As each of these relations can be obtained in an expected time of $O(\log(q)^{O(1)})$, we conclude that a tree with $N_{\max}$ edges can be constructed in an expected time of 
\[ O(q^{2 - \frac{2}{g}}) \; . \]

Note that the depth of the tree is always bounded by $s$. In particular, as $s = O(\log(q))$, the depth of the tree is also in $O(\log(q))$.

\subsubsection*{Construction of the matrix}

We now assume we have constructed a tree $T$ with $N_{\max} = \lceil q^{1 - \frac{1}{g} + \frac{1}{g^2}} \rceil$ edges.

Similarly to above let
\[
S := \big\{  P_1 + \cdots + P_g \in \Div^{g}(\mathcal{C}) | \; \forall i=1, \ldots, g : P_i \in \mathcal{F} \cup (T - \{ * \}) \big\} \; ,
\]
\[ S^{\ns} := S \cap D^{g, \ns} \; . \]
Then $S$ contains $\geq \frac{1}{g!} \cdot (\# \mathcal{F} + \# (T - \{* \}))^g \geq \frac{1}{g!} \cdot q^{g - 1 + \frac{1}{g}}$ elements. By the first two inequalities of (\ref{number-ns}), $S^{\ns}$ contains at least $\frac{3}{4g!} \cdot q^{g - 1 + \frac{1}{g}}$ elements. This means that the probability that a relation (\ref{relation}) splits into elements of the factor base or vertices of the tree is
\[ \geq \frac{3}{8 g!} \cdot q^{-(1 - \frac{1}{g})} \; . \]

The expected number of relations (\ref{relation}) which have to be considered until a matrix of size $\tilde{O}(\# \mathcal{F})$ is constructed is therefore in $O(q^{2 - \frac{2}{g}})$. As each relation can be computed in a time of $O(\log(q)^{O(1)})$ and the depth of the tree is in $O(\log(q)^{O(1)})$, this means that a matrix of size $\tilde{O}(\# \mathcal{F})$ can be constructed in a time of
\[ \tilde{O}(q^{2 - \frac{2}{g}}) \; . \]

\subsubsection*{Linear algebra}

The linear algebra takes place on a sparse matrix with $\tilde{O}(q^{1-1/g})$ columns and $\lceil q^{1 - 1/g} \rceil$ rows. (If the group order is square-free, a matrix with $\# \mathcal{F} +1= \lceil q^{1 - 1/g} \rceil +1$ rows suffices, but if the matrix is not square-free, according to the description in \cite{EG}, one constructs a larger matrix.)

Writing the relations in rows, as the tree has depth $O(\log(q))$, every row contains only $O(\log(q))$ non-zero entries.

We apply the algorithm in Section 4 of \cite{EG} to compute a vector $v$ over $\mathbb{Z}/\ell\mathbb{Z}$ with $v R = 0$. This algorithm now terminates in an expected time of
\[ \tilde{O}(q^{2 - \frac{2}{g}}) \; .\]

As argued in \cite[Section 4.5]{GTTD}, the double large prime variation does not affect the failure probability of the linear algebra computation, that is, the results of \cite[Section 4]{EG} still hold: After an expected number of $O(\log(q)^{O(1)})$ restarts of the construction of the matrix $R$, the linear algebra computation leads to the solution of the DLP.

\subsubsection*{Final result}

We have seen that the construction of the tree of large prime relations, the construction of the matrix $R$ and the application of the linear algebra algorithm in \cite[Section 4]{EG} all have a running time of $\tilde{O}(q^{2 - \frac{2}{g}})$. Moreover, we have argued that after an expected number of $O(\log(q)^{O(1)})$ restarts of the computation of the matrix $R$, the linear algebra computation leads the solution to the DLP. This means that the total running time is in 
\[ \tilde{O}(q^{2 - \frac{2}{g}}) \; , \]
in accordance with the theorem in the introduction.

\subsubsection*{Storage requirements}

Clearly there exists a function in $\tilde{O}(q^{1 - \frac{1}{g} + \frac{1}{g^2}})$ such that the storage requirements for the tree are bounded by this function for every run of the algorithm.

The storage requirements for the matrix are (for every run of the algorithm) bounded by a function in $\tilde{O}(q^{1 - \frac{1}{g}})$. Note again that this is the case because we restart the construction of the matrix every time the linear algebra computation fails instead of inserting a new row.

\bibliography{dlpv-literatur}

\begin{thebibliography}{10}

\bibitem{Di-sd}
C.\ Diem.
\newblock {An Index Calculus Algorithm for Plane Curves of Small Degree}.
\newblock In {\em Algorithmic Number Theory -- ANTS~VII}. Springer-Verlag,
  2006.
\newblock Forthcoming.

\bibitem{EG}
A.~Enge and P.~Gaudry.
\newblock {A general framework for subexponential discrete logarithm
  algorithms}.
\newblock {\em Acta.\ Arith.}, 102:83--103, 2002.

\bibitem{GTTD}
P.~Gaudry, E.~Thom\'{e}, N.~Th\'{e}riault, and C.~Diem.
\newblock {A double large prime variation for small genus hyperelliptic index
  calculus}.
\newblock accepted for publication in Math.~Comp., 2005.

\bibitem{Hess-subexp}
F.~He{\ss}.
\newblock Computing relations in divisor class groups of algebraic curves over
  finite fields.
\newblock Submitted, available under
  \textsf{http://www.math.tu-berlin.de/$\sim$hess}.

\bibitem{Hess-RR}
F.~He{\ss}.
\newblock {Computing Riemann-Roch spaces in algebraic function fields and
  related topics}.
\newblock {\em J.\ Symbolic Computation}, 11, 2001.

\bibitem{HI-RR}
M.-D. Huang and D.~Ierardi.
\newblock {Efficient Algorithms for the Riemann-Roch Problem and for Addition
  in the Jacobian of a Curve}.
\newblock {\em J.\ Symbolic Computation}, 18:519--539, 1994.

\bibitem{Mak2}
K.~Khuri-Makdisi.
\newblock Asymptotically fast group operations on jacobians of general curves.
\newblock available on the arXiv under math.NT/0409209, 2004.

\bibitem{Mak}
K.~Khuri-Makdisi.
\newblock Linear algebra algorithms for divisors on an algebraic curve.
\newblock {\em Math.\ Comp.}, 73:333--357, 2004.

\bibitem{LP}
H.W. Lenstra and C.~Pomerance.
\newblock A rigorous time bound for factorin integers.
\newblock {\em J. Amer. Math. Soc.}, 5, 1992.

\bibitem{Nag}
K.~Nagao.
\newblock {Improvement of Th\'eriault Algorithm of Index Calculus of Jacobian
  of Hyperelliptic Curves of Small Genus}.
\newblock Cryptology ePrint Archive, Report 2004/161,
  \textsf{http://eprint.iacr.org/2004/161}, 2004.

\bibitem{Pi-AV}
J.~Pila.
\newblock Frobenius maps of abelian varieties and fining roots of unity in
  finite fields.
\newblock {\em Math.\ Comp.}, 55:745--763, 1990.

\bibitem{Pi-Curves}
J.\ Pila.
\newblock Counting points on curves over families in polynomial time.
\newblock available on the arXiv under math.NT/0504570, 1991.

\bibitem{Schoof-EC}
R.\ Schoof.
\newblock Elliptic curves over finite fields and the compuation of square roots
  mod $p$.
\newblock {\em Math.\ Comp.}, 44:483--494, 1985.

\bibitem{St}
H.~Stichtenoth.
\newblock {\em {Algebraic Function Fields and Codes}}.
\newblock Springer-Verlag, Berlin, 1993.

\bibitem{Vo}
E.~Volcheck.
\newblock {Computing in the Jacobian of a Plane Algebraic Curve}.
\newblock In M.-D.~Huang L.~Adleman, editor, {\em Algorithmic Number Theory --
  ANTS~I}, volume 877 of {\em LNCS}, pages 28--40, Berlin, 1994.
  Springer-Verlag.

\end{thebibliography}

\bibliographystyle{plain}

\appendix

\section{On the number of special divisors}

The purpose of this appendix is to prove Proposition \ref{number-special-divisors} in the introduction.

Let $\mathcal{C}$ be a curve of genus $g$ over $\mathbb{F}_q$. 

Let $\Div^g(\mathcal{C})$ be the set of effective divisors of degree $g$ on $\mathcal{C}$, and let $D_0$ be an divisor of degree $g$ on~$\mathcal{C}$. We have the surjective map $\Div^g(\mathcal{C}) \longrightarrow \Cl^0(\mathcal{C}), D \mapsto [D] - [D_0]$. Note that the set of special divisors of degree $g$ is exactly the subset of $\Div^g(\mathcal{C})$ where the map to $\Cl^0(\mathcal{C})$ is not injective.

The number of special divisors is therefore bounded from above by $2 \, ( \# \Div^g(\mathcal{C}) - \# \Cl^0(\mathcal{C}))$, and it suffices to prove that $\# \Div^g(\mathcal{C}) - \# \Cl^0(\mathcal{C}) = O(q^{g-1})$.

\smallskip

We follow the exposition to the zeta-function in \cite{St}. 

Let $L= \prod_{i=1}^{2g} (1 - \alpha_i t) \in \mathbb{C}(t)$ be the $L$-polynomial of $\mathcal{C}$, let $A_n$ be the number of divisors of degree $n$, let $B_n$ be the number of prime divisors of degree $n$ on $\mathcal{C}$, and let
\[ S := \sum_{i=1}^{2g} \alpha_i \; .\]
As the $\alpha_i$ can be arranged such that $\alpha_i \alpha_{g +i} = q$ for all $i = 1, \ldots, g$, we have
\[ \# \Cl^0(\mathcal{C}) = L(1) = q^g - S \cdot q^{g-1} + O(q^{g-1}) \; .\]
We thus have to show that
\begin{equation}
\label{A_n-to-show}
A_g = q^g - S \cdot q^{g-1} + O(q^{g-1}) \; . 
\end{equation}
The fact that any divisor of degree $g$ can be expressed as a (up to permutations) unique sum of prime divisors implies
\[ A_g = \sum_{\underline{e}} \prod_r \binom{B_r + e_r -1}{e_r} \; , \]
where the sum runs over all $\underline{e} \in \mathbb{N}_0^{\,g}$ with $\sum_r e_r r = g$ and the products run over $r \in \{1, \ldots g \}$. We have
\[ B_1 = q + 1 - S \; ,\]
and by \cite[Proposition V.2.9]{St}, we have in particular
\[ B_r = \frac{1}{r} \cdot q^r + O(q^{r-1}) \]
for $r \geq 2$.

This implies that
\[ A_g = \sum_{\underline{e}} \frac{1}{e_1!} (q - S)^{e_1} \cdot \prod_{r \geq 2} \frac{1}{e_r!} \cdot \frac{1}{r^{e_r}} \cdot q^{r \cdot e_r} + O(q^{g-1}) \; ,\]
i.e.
\[ A_g = \sum_{\underline{e}} (\prod_{r} \frac{1}{e_r!} \cdot \frac{1}{r^{e_r}}) \cdot (q^g - e_1 \cdot S \cdot q^{g-1})+ O(q^{g-1}) \; .\]
In order to derive (\ref{A_n-to-show}) it remains to be shown that
\begin{equation}
\label{=1-to-show-Sg}
\sum_{\underline{e}} \prod_{r} \frac{1}{e_r!} \cdot \frac{1}{r^{e_r}} = 1
\end{equation}
and
\begin{equation}
\label{=1-to-show-S(g-1)}
\sum_{\underline{e}} e_1 \cdot \prod_{r} \frac{1}{e_r!} \cdot \frac{1}{r^{e_r}} = 1 \; .
\end{equation}

Equation (\ref{=1-to-show-Sg}) is equivalent to 
\begin{equation}
\label{to-show-Sg}
\sum_{\underline{e}} \prod_{r} \frac{g!}{e_r!} \cdot \frac{1}{r^{e_r}} = g! \; .
\end{equation}
This is true because for any $\underline{e} \in \mathbb{N}_0^{\, g}$ with $\sum_r e_r r = g$, the set of permutations in $S_g$ having exactly $e_r$ $r$-cycles has $\prod_{r} \frac{g!}{e_r!} \cdot \frac{1}{r^{e_r}}$ elements.

\smallskip
We come to Equation (\ref{=1-to-show-S(g-1)}). Note that we have a bijection
\[\begin{array}{c} 
\{ \underline{e} \in \mathbb{N}_0^{\, g} | \; \sum_r e_r r = g , \; e_1 \neq 0 \} \longrightarrow  \{ \underline{e}' \in \mathbb{N}_0^{\, g-1} | \; \sum_r e'_r r = g-1 \} \; , \\
\underline{e} \mapsto \underline{e}'
\end{array} \]
with $e'_1 = e_1-1$ and $e'_i = e_i$ for all $i = 1, \ldots, g-1$.

Equation (\ref{=1-to-show-S(g-1)}) is then equivalent to
\begin{equation}
\sum_{\underline{e}'} \prod_{r} \frac{(g-1)!}{e'_r!} \cdot \frac{1}{r^{e'_r}} = (g-1)! \; ,
\end{equation}
where the sum runs over all $\underline{e}' \in \mathbb{N}_0^{\, g-1}$ with $\sum_r e'_r r = g-1$ and the products run over $r \in \{1, \ldots, g-1 \}$. This is true by the same argument as the one for Equation (\ref{to-show-Sg}) (with $S_g$ substituted by $S_{g-1}$).

\bigskip

\bigskip

\bigskip

\noindent
Claus Diem, Universit\"at Leipzig, Fakult\"at f\"ur Mathematik und Informatik,
 Augustusplatz 10/11, 04109 Leipzig, Germany. diem@math.uni-leipzig.de

\end{document}